\documentclass[11pt]{article}%
\usepackage{amsmath}
\usepackage{amsfonts}
\usepackage{amssymb}
\usepackage{graphicx}%
\providecommand{\U}[1]{\protect\rule{.1in}{.1in}}

\begin{document}

\title{On a paper of Daskalopoulos and Sesum}
\author{Bennett Chow}
\date{}
\maketitle

\begin{quote}
{\scriptsize When you're ridin' sixteen hours and there's nothin' much to do.
And you don't feel much like ridin', you just wish the trip was through. From
`Turn the Page' by Bob Seger}
\end{quote}

This is an exposition of aspects of the result of Daskalopoulos and Sesum
\cite{DasSes} that any complete noncompact ancient solution $\left(
\mathcal{M}^{2},g\left(  t\right)  \right)  $ to Ricci flow with bounded $R>0$
and finite width must be the cigar soliton (we follow some of the main ideas
of their proof). Around the same time S.-C. Chu
\cite{Chu2dNoncompactAncientSurface} (based on earlier work of Shi and
Ni--Tam) proved that finite width follows from the former properties.

Suppose $f:\mathcal{M}\rightarrow\mathbb{R}$ satisfies $\Delta f=R$. By the
Bochner formula for $\Delta\!\left\vert \nabla f\right\vert ^{2}$ and
$\operatorname{Rc}=\frac{R}{2}g$, we have%
\[
\Delta\left(  R+\left\vert \nabla f\right\vert ^{2}\right)  =\frac{\left\vert
\nabla R+R\nabla f\right\vert ^{2}}{R}+2\left\vert \nabla\nabla f-\frac{1}%
{2}\Delta f\,g\right\vert ^{2}+\Delta R+R^{2}-\frac{\left\vert \nabla
R\right\vert ^{2}}{R}.
\]
For any smooth bounded domain $\Omega\subset\mathcal{M}$, by the trace Harnack
estimate and the divergence theorem%
\[
\int_{\Omega}\left(  \frac{\left\vert \nabla R+R\nabla f\right\vert ^{2}}%
{R}+2\left\vert \nabla\nabla f-\frac{1}{2}\Delta f\,g\right\vert ^{2}\right)
d\mu\leq\int_{\partial\Omega}\nu\left(  R+\left\vert \nabla f\right\vert
^{2}\right)  d\sigma.
\]

By Shi's local derivative estimate, $\left\vert \nabla R\right\vert \leq C$ on
$\mathcal{M}$. Since $\operatorname{inj}\left(  p\right)  \geq c>0$,
$\lim_{p\rightarrow\infty}R\!\left(  p\right)  \!=\!0$ at each time $t$ by the
Cohn-Vossen theorem or Hamilton's curvature bumps result. Shi's estimate again
yields $\lim_{p\rightarrow\infty}\left\vert \nabla R\right\vert \!\left(
p\right)  \!=\!0$. Fix $t$ and suppose the width of $g\!\left(  t\right)  $ is
finite. For any $p_{i}\rightarrow\infty$, there exist embeddings $\varphi_{i}$
such that $\left(  \varphi_{i}^{\ast}g\left(  t\right)  ,p_{i}\right)  $
subconverges in $C^{\infty}$ to a flat $\left(  \mathcal{S}^{1}\times
\mathbb{R},g_{\infty},p_{\infty}\right)  $ for if any limit is a flat
$\mathbb{R}^{2}$, then the width of $g\!\left(  t\right)  $ is infinite. Push
forward by $\varphi_{i}$ the geodesic circle containing $p_{\infty}%
\!\in\!\mathcal{S}^{1}\times\mathbb{R}$ to bound a disk $\Omega_{i}$ in
$\mathcal{M}$. Then $\cup_{i}\Omega_{i}=\mathcal{M}$ and $\partial\Omega
_{i}\rightarrow\infty$. By $\left\vert \partial\Omega_{i}\right\vert \leq C$,
$|\int_{\partial\Omega_{i}}\nu\!\left(  R\right)  d\sigma|\leq\left\vert
\partial\Omega_{i}\right\vert \sup_{\partial\Omega_{i}}\left\vert \nabla
R\right\vert \rightarrow0$ as $i\rightarrow\infty$. Assume $g=e^{-f}%
(dx^{2}+dy^{2})$. Then $\frac{\partial f}{\partial t}=\Delta f=R$. On
$\mathbb{R}^{2}-\left\{  0\right\}  $ let $g=vg_{c}$, where $g_{c}%
=\frac{dx^{2}+dy^{2}}{x^{2}+y^{2}}$ is isometric to $a^{-1}g_{\infty}$, $a>0$.
Let $p_{i}=\left(  \theta_{i},s_{i}\right)  $. By the finite width condition,
for each subcylinder $\mathcal{C}$ of length $1$ in $\mathcal{S}^{1}%
\times\lbrack s_{i}-\frac{\rho_{i}}{2},s_{i}+\frac{\rho_{i}}{2}]$ we have
$\inf_{x\in\mathcal{C}}v\left(  x\right)  \leq C_{1}.$ By Proposition
2.4\label{ALT PROOF?} in \cite{DasSes} there exists $c>0$ such that $v\geq c$
in $\mathcal{S}^{1}\times\left[  s_{i}-\rho_{i},s_{i}+\rho_{i}\right]  $,
where $\rho_{i}\rightarrow\infty$. Thus the Harnack inequality\label{TRUE?}
for almost harmonic functions implies that $v\leq C_{2}$ in $\mathcal{S}%
^{1}\times\lbrack s_{i}-\frac{\rho_{i}}{2},s_{i}+\frac{\rho_{i}}{2}].$ One has
uniform higher derivative bounds so that $(\mathcal{S}^{1}\times\lbrack
-\frac{\rho_{i}}{2},\frac{\rho_{i}}{2}],v\left(  \theta,s+s_{i}\right)
g_{c})$ subconverges\label{WHAT ABOUT A WAVE?} pointwise to a flat metric
$g_{\infty}^{\prime}=v_{\infty}g_{c}$ on $\mathcal{S}^{1}\times\mathbb{R}$
isometric to $ag_{c}.$ This implies that $v_{\infty}\equiv a.$ Let $\left(
r,\theta\right)  $ be polar coordinates and $s=\log r$. Then $g_{c}%
=ds^{2}+d\theta^{2}$ and $f=-\log v+2s$. Thus $\nabla\nabla f\rightarrow0$ on
$\partial\Omega_{i}$ and $\left\vert \nabla f\right\vert \leq3a^{-1/2}$ on
$\partial\Omega_{i}$, so that $\lim_{i\rightarrow\infty}\int_{\partial
\Omega_{i}}\nu(\left\vert \nabla f\right\vert ^{2})d\sigma=0$. Hence
$\int_{\mathcal{M}}(\frac{\left\vert \nabla R+R\nabla f\right\vert ^{2}}%
{R}+2\left\vert \nabla\nabla f-\frac{1}{2}\Delta f\,g\right\vert ^{2})d\mu=0$
and therefore $g\!\left(  t\right)  $ is a steady soliton, which by Hamilton's
classification result must be the cigar soliton.

\textbf{Acknowledgments.} Much obliged to Peng Lu, Jiaping Wang, and Bo Yang
for helpful discussions.

\end{document}